\documentclass[12pt]{article}
\usepackage{graphicx} 
\usepackage{amssymb}
\usepackage{amsmath}
\usepackage{tikz}
\usepackage{amsthm}
\usepackage{multicol}
\usepackage{mathrsfs}
\usepackage{hyperref}
\usepackage{array}
\usepackage{tikz-cd}

\numberwithin{equation}{section}

\newtheorem{theorem}{Theorem}

\newtheorem{proposition}[theorem]{Proposition}
\newtheorem{lemma}[theorem]{Lemma}
\newtheorem{definition}[theorem]{Definition}

\newtheorem{remark}[theorem]{Remark}
\newtheorem{example}{Example}

\newenvironment{customproof}[1]
    {\par\vspace{1em}\noindent\textbf{#1.}\;}
    {\hfill$\qed$\par\vspace{1em}}

\newcommand{\N}{\mathbb{N}}
\newcommand{\M}{\hat{M}}
\newcommand{\T}{\hat{T}}
\newcommand{\x}{\hat{x}}
\newcommand{\F}{\hat{F}}
\newcommand{\set}[2]{\left\{ #1 \; \left|\; #2  \right.\right\}}
\newcommand{\dnull}{\Delta_0}
\newcommand{\dhat}{\hat{\Delta}}
\newcommand{\R}{\widehat{R}}
\newcommand{\Ha}{\widehat{H}}
\newcommand{\Z}{\mathcal{Z}}
\newcommand{\G}{\mathcal{G}}
\newcommand{\U}{\mathcal{U}}

\title{Improved estimates of statistical properties in some non-uniformly hyperbolic dynamical systems}
\author{Péter Bálint\thanks{
Department of Stochastics, Institute of Mathematics,
Budapest University of Technology and Economics,
M{\H u}egyetem rkp. 3., H-1111, Budapest, Hungary; HUN-REN--BME Stochastics Research Group,
Budapest University of Technology and Economics,
M{\H u}egyetem rkp. 3., H-1111 Budapest, Hungary; Erd\H{o}s Center of the HUN-REN Alfr\'ed R\'enyi Institute of Mathematics, Re\'altanoda utca 13-15, H-1053, Budapest, Hungary;
\textit{balint.peter@ttk.bme.hu}. Support of the NKFIH Research Fund, grants 142169 and 144059, is thankfully acknowledged.},
Ábel Komálovics\thanks{
 Department of Stochastics, Institute of Mathematics,
Budapest University of Technology and Economics,
M{\H u}egyetem rkp. 3., H-1111, Budapest, Hungary; \textit{abelkomalovics@edu.bme.hu}. Project no. 3709 has been implemented with the support provided by the Ministry of Culture and Innovation of Hungary from the National Research, Development and Innovation Fund, financed under the EKÖP-24-3-BME-177 funding scheme.}}
\date{January 2025}

\begin{document}
\everymath{\displaystyle}
\maketitle

\begin{abstract}
    Building upon previous works by Young, Chernov-Zhang and Bruin-Melbourne-Terhesiu, we present a general scheme to improve bounds on the statistical properties (in particular, decay of correlations, and rates in the almost sure invariant principle) for a class of non-uniformly hyperbolic dynamical systems. Specifically, for systems with polynomial, yet summable mixing rates, our method removes logarithmic factors of earlier arguments, resulting in essentially optimal bounds. Applications include Wojtkowski's system of two falling balls, dispersing billiards with flat points and Bunimovich's flower-shaped billiard tables.
\end{abstract}

\section{Introduction}

In this paper we consider statistical properties for a class of non-uniformly hyperbolic dynamical systems. Our primary motivation is to study some billiard-like systems with polynomial, yet summable mixing rates (see Section~\ref{examples} for the specific examples) for which we manage to improve upon previous bounds on the rates of correlation decay and in the almost sure invariance. Nevertheless, we present a general scheme in the hope that it can be applied to many other models. This scheme reformulates and somewhat extends ideas of Chernov-Zhang (\cite{ChZh,improved estimates}) and Bruin-Melbourne-Terhesiu (\cite{Ian}) in such a way that it can handle the above mentioned polynomial, yet summable decay rates.
Let us comment briefly on some previous results – it is important to emphasize that we do not aim for a complete survey of the rich literature on this subject, just mention works that are directly related to our approach. In her seminal papers \cite{Young_edc,Young} Young introduced a tower method which has turned out to be particularly efficient for studying statistical properties of non-uniformly hyperbolic systems. When modeling a non-uniformly hyperbolic system with a Young tower, a tower map is constructed that is either isomorphic, or is an extension of the studied dynamical system. Also, there exist one-sided and two-sided versions of Young towers appropriate for modeling non-invertible and invertible systems, respectively. Nonetheless, in all cases, the return process to the base of the tower is of Markovian character, and the statistics of the return times imply upper bounds on the mixing rates of the system with respect to a natural invariant probability measure. \footnote{In the examples of Section~\ref{examples}, this measure is the unique ergodic absolutely continuous invariant probability measure, but the method is applicable to other situations, for example, SRB measures, as well .}

To apply Young towers for studying polynomial decay of correlations in some nonuniformly hyperbolic planar billiard examples, Markarian (\cite{Markarian}), Chernov and Zhang (\cite{ChZh}) designed a powerful scheme.
Following the terminology introduced in \cite{Ian}, we will call this scheme  Chernov-Markarian-Zhang (CMZ) structure. We refer to Appendix \ref{CMZ structures}, \ref{2sided CMZ} and \cite{Ian} for the precise definitions on these structures, just mention here the most important features. Given a billiard map or some other non-uniformly hyperbolic dynamical system $(M,T,\mu_M)$, a subset $\hat{M}\subset M$ is identified with  $\mu_M(\hat{M})>0$. It is important that typically, the first return process of $T$ to $\hat{M}$ (which will be referred from now on as the ``geometric return process’’) can be studied by explicit geometric methods, resulting in polynomial tails, but, in contrast to the return processes of Young towers, it is not Markovian. Yet, the first return map $\hat{T}:\hat{M}\to\hat{M}$ has strong hyperbolic properties and can be modelled by a Young tower with exponential tails. Accordingly, $(M,T,\mu_M)$ itself can be modeled by a Young tower with sub-exponential tails. Thus an important question is \textit{how to relate the statistics of the geometric return process to the statistics of the Young tower return process?} As discussed in (\cite{Markarian,ChZh}) even in the lack of any additional information, it is possible to obtain estimates on the Young tower return process, but this involves some logarithmic factors when compared to the geometric return process. These logarithmic factors then propagate to the upper bounds on the rate of correlation decay, which appears to be an artefact of the method in most situations. In \cite{improved estimates} Chernov and Zhang reconsidered some billiard examples -- in particular, stadia, semi-dispersing billiards derived from infinite horizon Lorentz gases and dispersing billiards with cusps -- and showed that, using some additional information on the geometric return process, the logarithmic factors can be eliminated. For all the systems studied in \cite{improved estimates}, the upper bound on the decay of correlations obtained this way is of the type $n^{-1}$. This rate is not summable, which  is closely related to the superdiffusive limit laws (non-standard $\sqrt{n \log n}$ scaling) arising in these systems (\cite{lower bounds,SV07,cusp_CLT}).

Our main goal in this paper is to \textit{reformulate and extend the results of \cite{improved estimates} to systems with polynomial, yet summable decay rates}. In particular, we consider Chernov-Markarian-Zhang structures with some additional properties, see (C1)-(C5) in Section~\ref{conditions},  and show that, under such conditions, the Young tower return process has the same type of asymptotics as the geometric return process. Consequently, for some examples – Wojtkowski’s system of two falling balls, dispersing billiards with flat points, and Bunimovich flowers, see Section~\ref{examples} for details – correlations are shown to decay at a rate $n^{-\alpha}$ for an appropriate  $\alpha>1$, which improves previously known bounds. Let us comment on some further features of our results.\\
\textit{Generality of the exposition.} As mentioned above, on top of the Chernov-Markarian-Zhang structure framework, some additional information are needed to show that the rates of the Young tower return process agree with the rates of the geometric return process. These information can be captured by the evolution properties of a class of invariant curves, and probability measures supported on these curves, which will be referred to as standard pairs. Standard pairs were introduced in the context of hyperbolic biliards by Chernov and Dolgopyat (see eg.~\cite[Section 7.4]{billiards}) and have been widely used to study their statistical properties. Here we include only the amount of information on standard pairs that is needed for our exposition to work -- we do not even explicitly assume full hyperbolicity of the dynamics – in the hope that our conditions, expressed in (C1)-(C5) of Section~\ref{conditions}, are flexible enough to be applied to many other situations.
We aim for generality also at the level of the decay rates treated. Indeed, the return time statistics does not have to be exactly of the form $n^{-\alpha}$, we can handle rates of the form $r(n)$ which is an arbitrary regularly varying function of index $-\alpha$ for some $\alpha>0$.\\
\textit{Exact correlation asymptotics.} Concerning non-uniformly hyperbolic systems with polynomial decay of correlations, an important problem is to establish, on top of the upper bounds of \cite{Young,ChZh,improved estimates}, lower bounds on the rates of mixing, and thus to identify the exact asymptotics of correlation decay. Significant results in this direction, using operator renewal theory, were obtained for Pommeau-Manneville type intermittent interval maps by Sarig \cite{Sa} and Gou\"{e}zel \cite{Gou}. The renewal theory approach was extended to Chernov-Markarian-Zhang structures by Bruin, Melbourne and Terhesiu in \cite{Ian}; which obtained  in particular exact correlation asymptotics  for planar billiard systems with superdiffusive behaviour (such as stadia or dispersing billiards with cusps). Nonetheless, the approach of \cite{Ian} exploited some information on the non-standard limit laws arising in these systems, and the case of polynomial, yet summable decay rates (still diffusive regime) remained open. In the present paper we fill this gap and prove exact correlation asymptotics for such systems, which includes in particular the examples of Section~\ref{examples}. Our results apply to the general class of dynamically H\"{o}lder continuous observables, as in \cite{Ian}, supported on the subset $\hat{M}$ for exact correlation asymptotics, yet with no restrictions on their support for upper bounds. Let us mention that \cite{optimal} also considered systems with non-summable decay of correlations (such as dispersing billiards with cusps) while here, as mentioned above, we focus on models with summable polynomial decay, for which no results on optimal rates have been stated earlier.\\
\textit{Rates in the almost sure invariance principle.} As mentioned above, in our main examples, although correlations decay polynomially, by summability of the rates, the Birkhoff sums obey a standard central limit theorem (CLT). Recently there has been much progress about the almost sure invariance principle (ASIP), a strengthening of the CLT stating that the process associated to the Birkhoff sums can be simultaneously realized and well approximated by Brownian motion, see \cite{ASIP} and references therein. In particular, in \cite{ASIP} it is proved that for systems modeled by Young towers, the rate of this approximation – that is, the rate of the ASIP – can be estimated in terms of the tails of the return times to the base of the tower. By better estimates on the return times, we also improve upon the rates of the ASIP specifically for the examples of Section~\ref{examples}, as conjectured in \cite[Remark 3.3]{ASIP}.\\
The rest of this paper is organized as follows. In Section~\ref{conditions}, the setting is introduced, which is that of Chernov-Markarian-Zhang structures with some additional properties expressed in the form of conditions (C1)-(C5).  In Section~\ref{results}, the main results of this paper are stated which relate the asymptotics of the Young tower return process to the  geometric return process under our conditions. In Section~\ref{corollaries},  some corollaries of the main results are derived which concern decay of correlations and rates in the ASIP. In Section~\ref{examples}, we include the discussion of the above mentioned specific examples as applications of our general results. In the Appendices, we recall some terminology on Young towers and CMZ structures, along with some properties of regularly varying functions.\\
Throughout this paper, we will use the following notations. For two functions $f,g: \mathbb{R}^+ \rightarrow \mathbb{R}$, $f(x) \ll g(x)$ holds if there is a $C>0$ and an $x_0\in\mathbb{R}$ such that $|f(x)| \leq C g(x) $ for all $x>x_0$.
Two functions are asymptotic, if $f(x) \ll g(x) \ll f(x)$. We will denote this by $f(x) \asymp g(x)$.

\section{Conditions}\label{conditions}

In this section, we list our assumptions about the system $T:M\rightarrow M$.
\\
\\
\textbf{(C1) Chernov-Markarian-Zhang structure}

We assume that the system has a two-sided Chernov-Markarian-Zhang structure, for which we use the definitions and the terminology in \cite{Ian}, and recall the details in the appendix.
Briefly, we have an invertible system $(M, T, \mu_M)$ and a subset $\M\subset M$, with $\mu_M(\M)>0$, such that the first return dynamics $(\M, \T, \mu_{\M})$ can be modeled by a two-sided Young tower $(\dhat, \F, \mu_{\dhat})$ with base $\dnull\subset \dhat$.
This induces a tower representation $\pi:(\Delta, F, \mu_{\Delta}) \rightarrow (M, T, \mu_M)$ over the same base $\dnull$.

We make the following additional assumptions. Although not included literally in the definition of two-sided CMZ structures
 in \cite{Ian}, these are assumed elsewhere in the exposition of \cite{Ian} or \cite{ASIP}, and satisfied in all our examples. We assume that
 \begin{itemize}
     \item the smaller tower $\dhat$ has exponential tails.
     \item the following two inequalities hold, which express that, between consecutive returns, distances of points in the tower $\Delta$ can be controlled by their respective distances in the base $\dnull$.
     \begin{itemize}
         \item  There exist
         constants $K > 0$, $\gamma_0 \in (0,1)$ such that
\begin{equation} \label{eq:DistIneqIan}
    d(T^\ell x, T^\ell y) \leq K \big( d(x, y) + \gamma_0^{s(x, y)} \big)
\end{equation}
for all $x, y \in \dnull$, $0 \leq \ell < h(x, T, \dnull)$.
\item For all $x, y \in \dnull$ with $s(x,y)\ge 1$  and any $0 \leq \ell < h(x, T, \dnull)$
     \begin{equation} \label{eq:DistIneqAlexey}
         d(T^\ell x, T^\ell y) \leq K \big( d(x, y) + d( F_0x, F_0 y) \big).
     \end{equation}
     \end{itemize}
 \end{itemize}
 Let us comment on these inequalities. The role of \eqref{eq:DistIneqIan} is solely to ensure that H\"older continuous observables are also dynamically H\"older continuous, see our Proposition~\ref{dynamically Hölder prop} or \cite[Proposition 7.3]{Ian}\footnote{In fact, it suffices if \eqref{eq:DistIneqIan} holds, instead of $d(x,y)$, for another metric $d_0(x,y)$ obtained as a power of $d$.}. The analogous inequality \eqref{eq:DistIneqAlexey} is assumed in \cite{ASIP} for similar reasons, see Remark \ref{r:whw} on their relation. We further note that \eqref{eq:DistIneqAlexey} is needed only for our ASIP result Theorem~\ref{ASIP theorem}.

 We will denote the level sets of the first return function $R$ of $\M$ with
\[
    R_k=\set{x\in\M}{R(x)=k+1}.
\]
\\
\textbf{(C2) Standard curves}

We assume that there is a family of curves $\mathcal{U}$ in $\M$ such that the $\T$-image of any $\gamma\in\mathcal{U}$ is a countable family of curves from $\mathcal{U}$. We call the elements of $\mathcal{U}$ standard curves.
\\
\\
\textbf{Standard pair}

We call $(\gamma, \mu)$ a standard pair, if $\gamma\in\mathcal{U}$ and $\mu$ is a probability measure supported on $\gamma$ which is absolutely continuous with respect to the arclength measure $m_{\gamma}$ on $\gamma$.
We assume that the $\T$-image of a standard pair $(\gamma, \mu_{\gamma})$ is a countable  family of standard pairs,\\
$\G=\set{\left(\gamma_j, \mu_j=(\T_*\mu|_{\gamma_j}(\gamma_j))^{-1}\T_*\mu|_{\gamma_j}\right)}{j\in J \subset \N}$, where
$\T_*\mu|_{\gamma_j}(\gamma_j)=\mu(\T^{-1}\gamma_j)$.
This way,
\[
\T_*\mu(H)=\sum_{j\in J} \mu_j(H)(\T_*\mu|_{\gamma_j}(\gamma_j)
\]
for any $H$ measurable meaning that there is a probability factor measure $\lambda_{\G}$  of $\T_*\mu$ on the index set $J$ with mass function $f_{\lambda_{\G}}(j)= \T_*\mu|_{\gamma_j}(\gamma_j)$.
\\
\\
\textbf{Standard family}

A $\G= \set{\left(\gamma_i, \mu_i\right)}{i\in I}$ is a standard family, if  $(\gamma_i, \mu_i)$ is a standard pair for all $i\in I$, and there is a probability factor measure $\lambda_{\G}$ on $I$.
This induces a probability measure $\mu$ supported on $\bigcup_i \gamma_i$ that satisfies
\[
    \mu(H)= \int_I \mu_i(H) \operatorname{d}\!\lambda_{\G}(i)
\]
for all $H$ measurable.
Then, the image of a standard family is again a standard family.
\\
\\
\textbf{(C3) Foliation of $R_k$}

We assume that the sets $R_k$  can be foliated with standard curves such that $\mu_{\M}|_{R_k}$ can be considered a standard family $\G_0^{(k)}$.
We would like to mention that this condition does not require $R_k$ to be connected.
\\
\\
\textbf{(C4) Bound on unstable widths}

We want to consider the unstable width of $R_k$. We say that the unstable width of $R_k$ is dominated by $f(k)$, if for any standard curve $\gamma \in \U$, $m_{\gamma}(\gamma \cap R_k) \ll f(k)$.
We only need an upper bound, thus domination is enough for us.

We will assume that there exists a $t>1$ such that  the unstable width of $R_k$ is dominated by $k^{-t}$.
\\
\\
\textbf{$\Z$ function}

For a standard family $\G=(\cup_{i\in I}\gamma_i, \mu)$ with probability factor measure $\lambda_{\G}$ let us define the function
\[
    r_i : \gamma_i \rightarrow \mathbb{R}^+, \quad r_i(x)=d_{\gamma_i}(x, \partial \gamma_i),
\]
which is just the length of the shorter segment of $\gamma_i$ cut by $x$.
Then, set
\[
    \Z(\G)=\sup_{\varepsilon>0}\frac{\mu(\cup_{i\in I}r_i^{-1}]0, \varepsilon[)}{\varepsilon}
\]
From our previous assumption, we can see that
\[
    \Z(\G_0^{(k)}) \gg k^{t}.
\]
\\
\textbf{(C5) Growth Lemma}

We assume that there is a $t\geq p>1$, such that for
$\G_0^{(k)}=(R_k, \mu|_{R_k})$, $\G_{m+1}^{(k)} = \T_* \G_m^{(k)}$,
\[
    \Z(\G_m^{(k)}) \ll k^{t-p}
\]
holds for all $k\ge 1$ and  $m \geq 1$.

If $\T$ has nice hyperbolic properties, then (C5) is satisfied. This follows from the observation that
if $\U$ is the set of unstable curves, we can relate (C5) to the classical version of the Growth Lemma, which is a well-known result for the class of dispersing billiards for example, see \cite[Section 5.9]{billiards}.
\begin{lemma}[Growth Lemma]\label{l:growthclassic}
    If $\G_0$ is a standard family, and $\G_{m+1}=\T\G_m$, then there exists a $C>0$ and a $\theta<1$ such that
    \[
        \Z(\G_{m})=C(\theta^m\Z(\G_0) + 1)
    \]
    for all $m \geq 1$.
\end{lemma}

If we have a first step expansion estimate
\[
    Z(\G_1^{(k)}) \ll k^{t-p},
\]
then (C5) easily follows from the classical version of the Growth Lemma.

\section{The main result \label{results}}

\subsection{Necessary definitions}

\textbf{Decay rates.} To present our scheme for studying subexponential decay of correlations, first we specify the class of decay rates that our method can handle.

As we will study upper and lower bounds on the mixing rates, we want the bounding functions to have some structure. We will use the following class of functions, which have a substantial literature, see \cite{regular variation} for example.
\begin{definition}
    A function $r:\mathbb{R}^+ \rightarrow \mathbb{R}^+$ is regularly varying with index $-\alpha\in\mathbb{R}$, if
    \[
        \lim_{x \rightarrow \infty}\frac{r(\lambda x)}{r(x)}=\lambda^{-\alpha}
    \]
    holds for all $\lambda>0$.
\end{definition}

\textbf{Notations on hitting and return times.}
We define the hitting time of $\dnull$ and the first return time of $\dhat$ as follows.

\begin{definition}
    For an $\x\in\dhat$, the hitting time of $\dnull$ is
    \[
        h (\x, \F, \dnull) = \min \set{n\in\N^+}{\F^{n}\x \in \dnull}.
    \]
    For $x\in\Delta$,
    \[
        h(x, F, \dnull) = \min \set{n\in\N^+}{F^{n}x \in \dnull}.
    \]
    Also,
    \[
        R(\x, F, \dhat) = R( \pi(\x)).
    \]
\end{definition}

With these functions, we can define the sets we need in order to state our main result.

\begin{definition}
    \[
        D_n = \set{\x \in \dhat}{h(\x, \F, \dnull) > n}
    \]
    \[
        H_n = \set{\x \in \dhat}{ R(\x, F, \dhat) > n}
    \]
    \[
        A_n = \set{x\in \Delta_0}{ h(x, F, \dnull) > n}
    \]
\end{definition}

Since $(\dhat, \F)$ has exponential tails by (C1), there exists a $\rho<1$ such that
\[
    \mu_{\dhat}(D_n) \ll \rho^n.
\]
We call $\mu_{\dhat}(H_n)$ the return tails of $\dhat$ or $\M$, and $\mu_{\dnull}(A_n)$ the first hitting time tails of $\Delta$ or $M$.
The decay of the measures of $A_n$ will determine the rate of the Young tower return process.

\subsection{Statement of the main result}

\begin{theorem}\label{main}
    If $T:M\rightarrow M$ satisfies (C1)-(C5), $r:\N^+ \rightarrow \mathbb{R}^+$ is regularly varying and $a > 0$, then
    \begin{itemize}
        \item [(a)]
                    $ \qquad
                        \mu_{\dhat}(H_n) \ll r(n)
                        \qquad
                        \text{ implies }
                        \qquad
                        \mu_{\dnull}(A_n) \ll r(n).
                    $
        \item [(b)]
                     $ \qquad
                        \mu_{\dhat}(H_n) \gg n^{-a} \;
                        \qquad
                        \text{ implies }
                        \qquad
                        \mu_{\dnull}(A_n) \gg \mu_{\dhat}(H_n).
                    $
        \item [(c)] $ \qquad
                        \mu_{\dhat}(H_n) \asymp r(n) \;
                        \qquad
                        \text{ implies }
                        \qquad
                        \mu_{\dnull}(A_n) \asymp r(n).
                    $
    \end{itemize}

\end{theorem}

We would like to remark that the exponent $a>0$ in $(b)$ is independent of the function $r(n)$ or it's index.

Part $(c)$ follows from part $(a)$ and part $(b)$ immediately.


\subsection{Some properties of regularly varying functions}

We work with this class of bounding functions, since it has been extensively studied, and it is a broad family that includes the usual return tail estimates encountered in the literature.
Some examples for regularly varying functions are
\begin{itemize}
    \item $x^{-\alpha}(\log x)^{\beta}$ for all $\alpha, \beta \in \mathbb{R}$,
    \item $ x^{-\alpha}\exp \left((\log x)^{\gamma}\right)$ for all $\gamma \in ]0, 1[$ and $\alpha\in \mathbb{R}$,
    \item $x^{-\alpha}\left(\log r(x)\right)^{\beta}$ for all $\alpha, \beta \in \mathbb{R}$ and regularly varying function $r$,
    \item $Cr(x)$ and $r(x)v(x)$ for all constants $C>0$ and regularly varying functions $r$ and $v$.
\end{itemize}
We will utilize two properties of regularly varying functions.

Firstly,
\begin{equation}\label{homogenous}
    r(\lambda x) \ll r(x)
\end{equation}
holds for any $\lambda>0$. This is easy to verify by definition.

\begin{lemma}\cite[Theorem 2.1.1]{regular variation}\label{index}
If $r:\mathbb{R}^+ \rightarrow \mathbb{R}^+$ is a regularly varying function with index $-\alpha$, then
\begin{equation*}
    \exists \lim_{x \rightarrow \infty} \frac{\log r(x)}{\log x} =  -\alpha.
\end{equation*}
\end{lemma}
We provide a proof in Appendix \ref{proof of rv}.
The second property we will use is a reformulation of the previous limit. Namely, that for all $\varepsilon>0$,
\begin{equation}\label{almost poly}
    x^{-\alpha-\varepsilon} \ll r(x) \ll x^{-\alpha+ \varepsilon}
\end{equation}
holds.

\subsection{Proof of Theorem \ref{main}}

For $x\in \Delta$, define the lap number as
\[
    L(x, \dhat) = \max \set{k\in\N}{\exists i_1, \dots, i_k\leq h(x, F, \dnull): \; F^{i_j}x\in \dhat \; \forall j}
\]
 and for a $b>0$, let
\[
    B_{n, b} = \set{x\in \Delta}{ L(x, \dhat) > b \log n}.
\]

\begin{customproof}{Proof of part $(a)$}
    We will divide $A_n$ into two disjoint subsets along $B_{n, b}$ and calculate the two measures separately.
    We showcase the bound on the measure of $A_n \cap B_{n, b}$ only for the sake of completeness, since it need not be improved from the one found in \cite{BBNV}, for example.
    If $x\in A_n \cap B_{n, b}\subset \dhat \cap B_{n,b}$, we can see that for at least n iterations, $x$ will stay outside of $\dnull$.
    Simultaneously, $x$ will hit $\dhat$ at least $b \log n$ times, meaning, that $\F^i(x) \notin \dnull$ for all $i\leq b \log n$.
    In other words,
    \[
        A_n \cap B_{n, b} \subset \dhat \cap B_{n,b} \subset D_{b \log n}.
    \]
    This implies
    \begin{align*}
        \mu_{\Delta}(A_n \cap B_{n, b})
        \ll&
        \mu_{\dhat}(D_{b \log n})
        \ll
        \theta^{b \log n}
        \ll r(n)
    \end{align*}
    for b large enough.

    As for $A_n \setminus B_{n, b}$, there is a natural combinatorial bound as well.
    We call consecutive $F$-images of a point $x$ avoiding $\dhat$ an excursion.
    We denote the longest excursion of $x$ in the first $h(x, F, \dnull)$ steps by $I(x)$. That is,
    \[
        I(x)=\max\set{R(F^kx)}{k< h(x, F, \dnull), \; F^kx\in\dhat}.
    \]
    If $x\in A_n \setminus B_{n, b}$, then $I(x)>\frac{h(x, F, \dnull)}{b\log n}>\frac{n}{b \log n}$. Thus, there exists a $k$ such that $\x=\T^kx$ is in $H_{\frac{n}{b \log n}}$.
    This means
    \[
        A_n \setminus B_{n, b}\subset \bigcup_{k \leq b \log n}\F^{-k} H_{\frac{n}{b \log n}}.
    \]
    Hence,
    \begin{align*}
        \mu_{\Delta}\left( A_n \setminus B_{n, b} \right)
        \ll&
        b \log n \, r\left(\frac{n}{\log n}\right)
        \ll
        \log n \left(\frac{n}{\log n}\right)^{-\alpha + \varepsilon}
        \\
        \ll&
        (\log n)^{1+\alpha - \varepsilon}n^{2\varepsilon} r(n)
        \ll
        n^{3\varepsilon}r(n)
    \end{align*}
    for all $\varepsilon>0$.

    It is clear we need a more careful analysis.
    We can see that the weak point of the above argument is the possibility of all the excursions of $x$ being the same order as the longest one.
    For a $q<1$, we define
    \begin{align*}
        \R_k &= \set{\x \in R_k}{\exists m\leq b \log k: \; \T^m \x \in \pi H_{k^q}} \subset \M;\\
        \Ha_k &=\bigcup_{m>k} \pi^{-1} \R_m\subset H_k\subset\dhat.
    \end{align*}

    \begin{proposition}\label{prop}
        Let $T:M\to M$ satisfy (C1)-(C5). Then, for any $b > 0$ there exists a $q < 1$ and a $\delta>0$ such that
        \begin{equation}\label{prop eq}
            \mu_{\M}(\R_k) \ll k^{-\delta} \mu_{\M}(R_k).
        \end{equation}

    \end{proposition}
     Heuristically, this means that for any point $\x$ the probability, that before $b \log n $ iterations of $\F$ we find another point the return time of which is comparable with that of $\x$, is small. Let us mention that relations analogous to \eqref{prop eq} are sometimes included among the assumptions of other works (see eg.~assumption (H3') in \cite{optimal}) while here we derive this property from our conditions (C1)-(C5).

    Recall, that for $x \in A_n\setminus B_{n, b}$, $I(x)> \frac{n}{b \log n}=n'$ denotes the length of the longest excursion of $x$ before time $n$.
    We separate two main cases.

        If $I(x)\leq \frac{n}{10}$, then the second longest excursion must be at least $\frac{9}{10}\frac{h(x, F, \dnull)}{b \log n}>\frac{9}{10}\frac{n}{b \log n}= \frac{9}{10}n'>(n')^q$ long for any $q<1$ and $n$ large enough, meaning
    \[
        \set{x \in A_n\setminus B_{n, b}}{I(x)\leq\frac{n}{10}}
        \subset
        \Ha_{n'}.
    \]

    Hence, by Proposition \ref{prop},
    \begin{align*}
        \mu_{\dhat}\set{x \in A_n\setminus B_{n, b}}{I(x)\leq\frac{n}{10}}
        \ll&
        (n')^{-\delta} r(n')
        =
        \left(\frac{n}{b \log n}\right)^{-\delta} r\left(\frac{n}{b \log n}\right)
        \\
        \ll&
        \left(\frac{n}{b \log n}\right)^{-\delta}\left(\frac{n}{b \log n}\right)^{-\alpha +\varepsilon}
        \\
        \ll&
        n^{-\delta+2\varepsilon}(\log n)^{\delta + \alpha -\varepsilon} r(n)
    \end{align*}
    holds for all $\varepsilon > 0$. If we choose $\varepsilon < \delta/2$, we get
    \[
        \mu_{\dhat}\set{x \in A_n\setminus B_{n, b}}{I(x)\leq\frac{n}{10}}
        \ll
        r(n).
    \]

    If $I(x)>\frac{n}{10}$, then we will use the following argument from the proof of \cite[Proposition 9.7]{decay of flows}.
    For any $(x, 0)\in \dnull$ let $(x, l_1(x))\in\dhat$ be the starting point of the longest excursion, meaning $I(x, 0)=R(x, l_1(x))$.
    Then,
    \begin{align*}
        \mu_{\dnull}&\set{p \in A_n}{I(p)>\frac{n}{10}}
        \\
        \asymp&
        \mu_{\dhat}\set{(x, 0) \in A_n}{R(\T^{l_1(x)}x)>\frac{n}{10}}
        \\
        =&
        \mu_{\dhat}\set{(x, l_1(x)) \in \dhat}{R(\T^{l_1(x)}x)>\frac{n}{10}}
        \\
        =&
        \mu_{\dhat}\set{(x, l_1(x)) \in \dhat}{R\circ\pi(x, l_1(x))>\frac{n}{10}}
        \\
        \leq&
        \mu_{\dhat}\set{p \in \dhat}{R\circ\pi (p)>\frac{n}{10}}
        \\
        \leq&
        \mu_{\dhat}(H_{\frac{n}{10}})
        \ll
        r(n/10)
        \ll
        r(n).
    \end{align*}

Hence, to complete the proof of part (a) of Theorem~\ref{main}, it remains to prove Proposition~\ref{prop}.

\begin{customproof}{Proof of Proposition \ref{prop}}
    We want to prove
    \[
        \frac{\mu_{\M}(\R_k)}{\mu_{\M}\left(R_k\right)}
        \ll
        k^{-\delta},
    \]
    for some $\delta>0$ where
    \[
        \R_k = \set{\x \in R_k}{\exists m\leq b \log k: \; \T^m \x \in \pi H_{k^q}}
        =
        \bigcup_{m\leq b \log k}R_k\cap \T^{-m} \pi H_{k^q}.
    \]
    It is enough to show that
    \begin{equation} \label{single}
        \mu_{\M}\left(R_k\cap \T^{-m} \pi H_{k^q}\right)
        \ll
        k^{-\delta'}\mu_{\M}\left(R_k\right)
    \end{equation}
    holds with some $\delta'>0$ for all $m\leq b \log k$. Indeed, summation on $m$ then gives the statement of the Proposition for any $\delta>0$ with $\delta'>\delta>0$, as the number of terms is $b \log k \ll k^{\delta'-\delta}$ for all such $\delta$.

    To prove \eqref{single}, we fix some $m\leq b \log k$ and compute
    \begin{align}\label{conditional}
        \frac{\mu_{\M}\left(R_k\cap \T^{-m} \pi H_{k^q}\right)} {\mu_{\M}(R_k)}
        =&        \frac{\mu_{\M}\left(\T^{-m}\left(\T^m R_k\cap \pi H_{k^q}\right)\right)} {\mu_{\M}(\T^{-m}\T^m R_k)} \nonumber
        =
        \mu_{\M}\left(\pi H_{k^q} \;|\; \T^m R_k\right)
        \\
        =&
        \mu_{\M}\Big(\bigcup_{l\geq k^q} R_l \;\left|\; \T^{m} R_k\right.\Big).
    \end{align}

    We may regard the conditional measure $\mu_{\M}|_{R_k}$ and its push-forward by $\T^m$,
    $\mu_{\M}|_{\T^mR_k}$ as standard families $\G^{(k)}_0$ and $\G^{(k)}_m$, respectively. Hence, we can estimate the conditional probability in \eqref{conditional} by using the unstable curves $\gamma$ included in the standard family $\G^{(k)}_m$.
    The measures $\mu_{\gamma}$ on these unstable curves are uniformly equivalent to the Lebesgue (arclength) measure $m_{\gamma}$ on $\gamma$, for each standard pair $(\gamma, \mu_{\gamma})$. Hence, we only have to compare the Lebesgue measures of $\gamma\cap (R_l \cap \T^mR_k)$ and $\gamma\cap \T^mR_k$.
    Then, integration with respect to the factor measure $\lambda=\lambda_{\G^{(k)}_m}$ will provide an estimate of (\ref{conditional}).
    By our assumption (C4), the unstable size of $R_l$ is uniformly bounded from above, which gives $\gamma\cap (R_l \cap \T^mR_k) \ll l^{-t}$. Thus, for any $\gamma_i \in \G^{(k)}_m$
    \[
        \frac{m_{\gamma_i}(\gamma_i\cap (R_l \cap \T^mR_k))}{m_{\gamma_i}(\T^mR_k)}\ll \frac{l^{-t}}{m_{\gamma_i}(\T^mR_k)}.
    \]
    Hence,
    \begin{align*}
        \mu_{\M}\Big(\bigcup_{l\geq k^q} R_l \;\left|\; \T^m R_k\right.\Big)
        =&
        \sum_{l\geq k^q} \left(\int_I \frac{m_{\gamma_i}(\gamma_i\cap (R_l \cap \T^mR_k))}{m_{\gamma_i}(\T^mR_k)} \operatorname{d}\!\lambda(i) \right)
        \\
        \ll&
        \sum_{l\geq k^q} \left(
        \, l^{-t} \int_I (m_{\gamma_i}(\T^mR_k))^{-1} \operatorname{d}\!\lambda(i) \right)
        \\
        \ll&
        \,  \sum_{l\geq k^q} \left( l^{-t}\Z(\G^{(k)}_m) \right)
        \ll
        \frac{\sum_{l\geq k^q}l^{-t}}{k^{p-t}}
        \ll
        \frac{k^{q(1-t)}}{k^{p-t}}
        \\
        =&
        k^{(q-1)(1-t)+1-p}.
    \end{align*}

    Here $q-1$, $1-p$ and $1-t$ are all negative. Yet, by choosing $q$ sufficiently close to $1$, in particular for $1-\frac{p-1}{t-1}<q<1$, the above computation gives an estimate $\ll k^{-\delta'}$  on \eqref{conditional}, thus providing \eqref{single} and completing the proof of Proposition~\ref{prop}.
    \end{customproof}\end{customproof}

\begin{customproof}{Proof of part $(b)$}
    Let us fix a $b>0$ satisfying $b\log \theta <-a$. We would
    like to compare the measure of $H_n\subset\dhat$ to the measure
    of $A_n\subset \dnull$. For $y\in\dnull$, the column above $y$
    is the set of points $\F^ky, k=0,\dots, h(y,\F,\dnull)-1$.
    Also, let $\pi_{\dnull}:\dhat\to\dnull$ denote the projection
    to the base of the tower, which maps the column of $y$ to $y$.

    By Proposition \ref{prop} there exists some $C>0$ such that,
    \begin{align}
        \mu_{\dhat}(H_n)
        \ll&
        (1-C n^{-\delta})\mu_{\dhat}(H_n)
        \ll
        \mu_{\dhat}(H_n\setminus\Ha_n) \nonumber
        \\
        \ll&
        \mu_{\dhat}\left((H_n \setminus\Ha_n)\cap B_{n, b}\right)
        +
        \mu_{\dhat}\left((H_n \setminus\Ha_n)\setminus B_{n, b}\right)
    \end{align}

    As for these two measures,
    \begin{align*}
        \mu_{\dhat}\left((H_n \setminus\Ha_n)\cap B_{n, b}\right)
        \ll&
        \mu_{\dhat}(\dhat \cap B_{n, b})
        \ll
        \mu_{\dhat}(D_{b\log n})
        \ll
        \theta^{b\log n}
        =
        n^{b\log\theta}.
    \end{align*}

    \begin{align*}
        &\mu_{\dhat}\left((H_n \setminus\Ha_n)\setminus B_{n, b}\right)
        =\\
        &=\mu_{\dhat}\set{x \in H_n}{h(x, \F, \dnull)\leq b\log n , \, \nexists m\leq b\log n: \, \F^m x \in H_{n^q}}
    \end{align*}

    For any such $x\in (H_n \setminus\Ha_n)\setminus B_{n, b}$, let $y=\pi_{\dnull}(x)\in \dnull$ denote its projection to $\dnull$, the base of the tower.
    From the right hand side, it is clear that $x$ is the uppermost $F$-image of $y$ in it's column satisfying $x\in H_n$ and $h(x, \F, \dnull)\leq b\log n$ simultaneously. Also, for any $z \in \dnull$ there exists at most one point in the column of $z$ that belongs to $(H_n \setminus\Ha_n)\setminus B_{n, b}$. Let $E_n\subset\dnull$ be the set of all points $y\in\dnull$ for which there exists (exactly one) $x\in (H_n \setminus\Ha_n)\setminus B_{n, b}$. Then, $\pi_0$ is a bijection between $(H_n \setminus\Ha_n)\setminus B_{n, b}$ and $E_n$. Hence,

    \begin{align*}
        \mu_{\dhat}\left((H_n \setminus\Ha_n)\setminus B_{n, b}\right)
        =&
        \mu_{\dhat}(\pi_0^{-1}(E_n))
        =
        \mu_{\dhat}(E_n)
        \\
        \ll&
        \mu_{\dhat}\set{y\in \dnull}{\exists  l(y) < h(y, \F, \dnull): \; x = \F^{l(y)}y\in H_n }
        \\
        \ll&
        \mu_{\dhat}\set{y\in \dnull}{I(y)>n}
        \\
        \ll&
        \mu_{\dhat}(A_n)
        \ll
        \mu_{\dnull}(A_n).
    \end{align*}

    Combining these three upper bounds, we can see that
    \begin{align*}
        \mu_{\dhat}(H_n) - n^{b\log \theta}
        \ll&
        \mu_{\dhat}(H_n) - \mu_{\dhat}\left((H_n \setminus \Ha_n ) \cap B_{n, b}\right)
        \ll
        \mu_{\dhat}\left((H_n \setminus \Ha_n ) \setminus B_{n, b}\right)
        \\
        \ll&
        \mu_{\dnull}(A_n).
    \end{align*}

Since $b\log \theta <-a$, this implies
\[
    \mu_{\dhat}(H_n) \ll  \mu_{\dnull}(A_n).
\]
\end{customproof}

\section{Corollaries and applications \label{corollaries}}

\subsection{Upper bound on the decay of correlation for general obervables}\label{decay of correlation}
In Sections \ref{decay of correlation} and \ref{asymptotic decay} we take a look at the decay of correlations of observables and ways of utilizing Theorem \ref{main} to estimate it.
With this in mind, let us introduce the notations of the following sections.

For functions $f, g\in L^1(M)$, let
    \[
        C_T(f, g, n)=\int_M f \cdot (g\circ T^n) d\mu - \int_M f d\mu \cdot \int_M g d\mu
    \]
    denote the correlation of $f$ and $g$ at time $n$.

It is well konwn that $\mu_{\dnull}(A_n)$ gives an upper bound on the rate of decay correlations in $n$ for sufficiently regular observables. Let us recall the standard notion of H\"older continuity.

\begin{definition}
    Let $(M, d)$ be a bounded metric space equipped with a Borel probability measure $\mu$. Fix $\eta\in ]0, 1[$.
    For a function $f\in L^{\infty}(M)$  define
    \[
    || f ||_{\eta}= |f|_{\infty} + |f|_{\eta}, \quad
    |f|_{\eta}=\sup_{\substack{x, y \in M \\ x \neq y}}\frac{|f(x)-f(y)|}{d(x, y)^{\eta}}.
    \]
    $f$ is H\"older continuous with exponent $\eta$ if its H\"older norm, $|| f ||_{\eta}$ is finite.
\end{definition}

As in \cite{Ian}, we consider the (more general) class of dynamically H\"older continuous functions.  Recall from Appendix \ref{2sided CMZ}, or \cite{Ian} that if $T: M \to M$ has a two-sided Chernov-Markarian-Zhang structure, then there is a separation time function $s$ on $\dnull$. In the one-sided case it is separating, but in the two-sided case, given a non-trivial factor map $\pi'$, it's lift does not separate the points of the $\dnull$.

\begin{definition}
    Suppose $T:M \rightarrow M$ admits a two-sided Chernov-Markarian-Zhang structure, as in (C1). Fix $\vartheta<1$. For a function $f:M \rightarrow \mathbb{R}$ define
    \[
        ||f||_\mathcal{H}=|f|_{\infty} + |f|_{\mathcal{H}},
        \quad
        |f|_{\mathcal{H}}
        =
        \sup_{\substack{x, y\in\dnull \\ x\neq y}}
        \sup_{0\leq k < h(x)}
        \frac{|f(F^k x) - f(F^k y)|}{d(x, y) + \vartheta^{s(x, y)}}.
    \]
    We say that f is dynamically Hölder, if $||f||_\mathcal{H}<\infty$ and denote by $\mathcal{H}(M)$ the space of such observables.
    Write $\mathcal{H}(\M)=\set{f\in\mathcal{H}(M)}{ \operatorname{supp}f\subset \M}$.
\end{definition}

    The authors of \cite{Ian} give sufficient conditions for $T:M \rightarrow M$ to satisfy in order for the Hölder observables to be dynamically Hölder as well.

\begin{proposition}\cite[Proposition 7.3]{Ian}\label{dynamically Hölder prop}
    Let $\eta \in (0,1]$ and let $d_0$ be a bounded metric on $M$. Let $C^\eta(M)$ be the space of observables that are $\eta$-H\"older with respect to $d_0$. Suppose that there exist $K > 0$, $\gamma_0 \in (0,1)$ such that
\[
d_0(T^\ell x, T^\ell y) \leq K \big( d_0(x, y) + \gamma_0^{s(x, y)} \big)
\]
for all $x, y \in \dnull$, $0 \leq \ell < h(x, T, \dnull)$.

Then $C^\eta(M) \subset \mathcal{H}(M)$ where we may choose any $\theta \in [\gamma_0^\eta, 1)$ and $d = d_0^{\eta'}$ for any $\eta' \in (0,\eta]$.

\end{proposition}
Hence, by inequality~\eqref{eq:DistIneqIan} assumed in condition (C1), Hölder observables with respect to the metric $d(x,y)$ are also dynamically Hölder for all systems satisfying (C1).

Young proved in \cite[Theorem 3]{Young} that if a one-sided Young tower $(\Delta', F', \mu_{\Delta'})$ satisfies $\mu_{\dnull}(A_n) \ll n^{-\alpha}$ for some $\alpha > 1$, then
\[
    |C_F(f, g, n)| \ll \sum_{k\geq n}k^{-\alpha} (\ll n^{-\alpha+1})
\]
holds for all Hölder observables $f,g\in C^\eta(\Delta')$.

\cite[Theorem 2.10]{coupling argument} proves that the previous bound also holds for two-sided Young towers and Hölder observables. As noted in \cite[Section 7.2]{Ian},
the proof works for dynamically Hölder observables as well. Indeed, the only statement that relies on the properties of the observables is \cite[Proposition 5.3]{coupling argument}, which holds in the dynamically Hölder case as well.
Moreover, here we argue that the proof also works if, instead an exact polynomial upper bound $n^{-\beta}$ rate, we have an arbitrary regularly varying function $r(n)$ with with index $-\alpha<-1$, see the statement of Theorem~\ref{correlation general} below.
Indeed, the only two places where the authors rely on the attributes of the bounding function $n^{-\beta}$ are the proofs of \cite[Theorem 2.10]{coupling argument} and \cite[Proposition 4.10]{coupling argument} (on which the proof of \cite[Theorem 2.10]{coupling argument} relies).
There, the only relevant property of polynomial functions - besides summability - is the fact that there is a $\gamma\in\mathbb{R}$ such that
\[
    \sum_{m\geq n}r\left(\frac{m}{k}\right) \leq k^{\gamma}\sum_{m\geq n}r(m) \quad \forall k\in\N^+
\]
for all $n$ large enough.
This holds for any regularly varying function $r$ with index $-\alpha<-1$, since it is true for every element of the sum.

Hence, we get that \cite[Theorem 2.10]{coupling argument} also holds for dynamically Hölder observables and regularly varying bounding functions with index $-\alpha<-1$. Theorem~\ref{correlation general} will also use the following simple Lemma (proof included in Appendix \ref{proof of rv}).
\begin{lemma}\label{asymptotic}
    If $r:\mathbb{R}^+ \rightarrow \mathbb{R}^+$ is a non-increasing regularly varying function with parameter $-\alpha<-1$, then
    \[
    \sum_{k\geq n} r(k) \asymp \int_n^{\infty}r(x)\operatorname{d}x.
    \]
\end{lemma}
Altogether, we get the following result.
\begin{theorem}\label{correlation general}
    Let $T: M\rightarrow M$ satisfy (C1)-(C5). If
    \[
        \mu_{\dhat}(H_n)\ll r(n)
    \]
    holds for some regularly varying function $r$ with index $-\alpha<-1$, then
    \[
        |C_T(f, g, n)| \ll \sum_{k\geq n} r(k) \asymp \int_n^{\infty}r(x)\operatorname{d}x
    \]
    for all $f, g \in \mathcal{H}(M)$.
\end{theorem}

\subsection{Asymptotic bounds on the decay of correlation for special observables}\label{asymptotic decay}

If we restrict the support of the observables to the exponentially mixing $\M$, we can get asymptotic bounds on the correlations as well.

Recall from Appendix \ref{2sided CMZ}, that
\[
    \overline{h}=\int_{\dnull}h\operatorname{d}\!\mu_{\dnull}.
\]

\begin{theorem}\label{correlation}
    Let $T: M\rightarrow M$ satisfy (C1)-(C5). If
    \[
        \mu_{\M}(H_n)\asymp r(n)
    \]
    holds for some regularly varying function $r$ with index $-\alpha<-1$, then
    \begin{itemize}
        \item [(a)]
                \[
                    |C_T(f, g, n)|\asymp
                    \overline{h}\sum_{k\geq n}\mu_{\dnull}(A_k)
                    \left| \int_M f \operatorname{d}\! \mu\int_M g \operatorname{d}\!\mu \right|
                    \asymp
                    \sum_{k\geq n}r(k)
                    \asymp
                    \int_n^{\infty}r(x) \operatorname{d} x
                \]
                for all $f, g\in \mathcal{H}(\M)$ with $\int_M f \operatorname{d}\! \mu\int_M g \operatorname{d}\!\mu \neq 0$.
        \item [(b)] For any $\varepsilon>0$,
                \[
                    |C_T(f, g, n)|
                    \ll
                     ||f||_{\mathcal{H}}||g||_{\mathcal{H}}n^{-\alpha + \varepsilon}\log n
                \]
                holds for all $f, g\in \mathcal{H}(\M)$ with $\int_M f \operatorname{d}\! \mu=0$.
    \end{itemize}
\end{theorem}
We included the $\log n$ factor in $(b)$, since if we have a function $r$ such that (\ref{almost poly}) holds with $-\alpha$ instead of $-\alpha+\varepsilon$, then actually,
\[
    |C_T(f, g, n)|
    \ll
    ||f||_{\mathcal{H}}||g||_{\mathcal{H}}n^{-\alpha}\log n
\]
holds true.
\begin{customproof}{Proof}
    By Theorem \ref{main}, $\mu_{\dnull}(A_n)\asymp r(n)$ holds.
    Since $f$ has an index $-\alpha<-1$, we know that for all $\varepsilon>0$, $\mu_{\dnull}(A_n)\asymp r(n) \ll n^{-\alpha + \varepsilon}$.
    By \cite[Theorem 7.4 $(a)$]{Ian} and \cite[Proposition 3.2]{Ian},
    \begin{equation}\label{corr bound}
        \left|C_T(f, g, n)-\overline{h}\sum_{k\geq n}\mu_{\dnull}(A_k)\int_M f \operatorname{d}\! \mu\int_M g \operatorname{d}\!\mu \right|
        \ll
        ||f||_{\mathcal{H}}||g||_{\mathcal{H}}(n^{-\alpha + \varepsilon}\log n + \zeta_{\alpha -\varepsilon}(n) )
    \end{equation}
    holds for all $f, g\in \mathcal{H}(\M)$, where
    \[
        \zeta_{a}(n)
        =
        \begin{cases}
            n^{-a} & a > 2 \\
            n^{-2} \log n & a = 2 \\
            n^{-2(a-1)} & 1 < a < 2
        \end{cases}.
    \]
    From Lemma \ref{asymptotic}, it is clear that
    \[
        \sum_{k\geq n}\mu_{\dnull}(A_k)
        \asymp
        \sum_{k\geq n}r(k)
        \asymp
        \int_n^{\infty}r(x) \operatorname{d} x
        \gg
        n^{-\alpha -\varepsilon + 1}
    \]
    Thus, in order to prove our statement, it is sufficient to check that both terms on the right side of (\ref{corr bound}) are decaying faster than $n^{-\alpha -\varepsilon + 1}$. Let us take $\varepsilon < 1/2$. Then,
    apparently
    \[
        n^{-\alpha +\varepsilon}\log n
        \ll
        n^{-\alpha -\varepsilon+1},
    \]
    thus we only have to take a closer look at $\zeta_a$. It is easy to check that
    \[
        \zeta_{\alpha -\varepsilon}(n) \ll
        n^{-\alpha -\varepsilon + 1}
    \]
    holds for all $\alpha - \varepsilon > 1$. This proves $(a)$.
    If $\int_M f \operatorname{d}\!\mu_M=0$, then the statement follows from \cite[Theorem 7.4 $(b)$]{Ian}, \cite[Proposition 3.2]{Ian} and part $(a)$ of Theorem \ref{main}.
\end{customproof}
If one does not have an asymptotic, only an upper bound on the tails of the return time to $\M$, then we can see that part $(a)$ of Theorem \ref{correlation}  offers no improvement over Theorem \ref{correlation general}, but regarding centered observables, it is still a better estimate.
\begin{remark}
    Let $T: M\rightarrow M$ satisfy (C1)-(C5). If
    \[
        \mu_{\M}(H_n)\ll r(n)
    \]
    holds for some regularly varying function $r$ with index $-\alpha<-1$, then,
      for any $\varepsilon>0$,
                \[
                    |C_T(f, g, n)|
                    \ll
                     n^{-\alpha +\varepsilon}\log n
                \]
                holds for all $f, g\in \mathcal{H}(\M)$ with $\int_M f \operatorname{d}\! \mu=0$.
\end{remark}

As in Theorem \ref{correlation} $(b)$, we include the logarithmic factor, since if $r(n) \ll n^{-\alpha}$, we can discard the $\varepsilon$ in the exponent.

\subsection{Almost sure invariance principle}\label{ASIP}

     We recall the following definition (see eg.~\cite{ASIP}).
    \begin{definition}
    Let $\kappa:\N^+\to \mathbb{R}$ denote a function with $\kappa(n)=o(\sqrt{n})$, and $c^2$ some positive constant. A random process $(S_n)_{n\ge 0}$ satisfies the almost sure invariance principle (ASIP) with variance $c^2$ and rate $o(\kappa(n))$ if $(S_n)_{n\ge 0}$ can be redefined on some probability space supporting a Brownian motion $W_n$ with variance $c^2$ such that
    \[
    S_n=W_n + o(\kappa(n)) \qquad \text{almost surely.}
    \]
    \end{definition}

    In \cite{ASIP} the authors study the almost sure invariance principle for a class of non-uniformly hyperbolic maps. Specifically, if our condition (C1) holds, then the system also satisfies the conditions formulated \cite[Section 2]{ASIP}\footnote{Let us recall that inequality \eqref{eq:DistIneqAlexey} included in condition (C1) is used only in Section~\ref{ASIP}}.

    The authors have a general result, \cite[Theorem 3.1]{ASIP} on the rates of the ASIP.
    As part $(a)$ of their theorem requires a bound on $\mu_{\dnull}(A_n)$, we can make an immediate improvement on it. (We note that the return time $h:\dnull\to\N^+$ is denoted by $\tau$ in \cite{ASIP}.)
     \begin{theorem}\label{ASIP theorem}
        Let $T : M \rightarrow M$ satisfy (C1)-(C5).
        Let $\varphi : M \rightarrow \mathbb{R}$ be Hölder continuous and centered, i.e. $\int \varphi \, d\mu_M = 0$.
        Let $S_n(\varphi) = \sum_{k=0}^{n-1} \varphi \circ T^k$.
        Assume that $\int h^2 \, d\mu_{\dnull} < \infty$.
        Then
        \begin{itemize}
            \item the limit $c^2 = \lim_{n \rightarrow \infty} n^{-1} \int |S_n(\varphi)|^2 d\mu_{\dnull}$ exists and
            \item if $\mu_{\dhat}(H_n) \ll n^{-\beta} (\log n)^{\gamma}$, with $\beta > 2$ and $\gamma \in \mathbb{R}$, then:\\ for any $\varepsilon > 0$ the process $S_n$ satisfies the ASIP with variance $c^2$ and rate $o(n^{1/\beta} (\log n)^{(1+\gamma)/\beta + \varepsilon})$.
        \end{itemize}

    \end{theorem}

\section{Examples}\label{examples}
In this section, we present three examples of systems for which our results improve upon earlier estimates. In particular, previous bounds on the rates of the decay of correlations were suboptimal for all three examples, and we also provide improved ASIP rates.

\begin{example}\textbf{Wojtkowski's system of two falling balls}\label{fallball}

    In Wojtkowski's system we have two balls moving on a vertical halfline. The end of the halfline is considered to be the floor. The two balls are accelerated by gravity and collide with each other and the floor elastically. If $m_1$ and $m_2$ denote the masses of the lower and the upper ball, respectively, then in case $m_1>m_2$ (when the lower ball is heavier) this system is hyperbolic and ergodic with respect to the Liouville measure (when restricted to the surface of constant total energy). Equivalently, when considering the Poincar\'e section of all elastic collisions (both between the two balls and between the lower ball and the floor), it is possible to introduce natural coordinates for which the Poincar\'e map, $T:M\to M$  is area preserving, non-uniformly hyperbolic and ergodic. We refer to \cite{BBNV} for further details.

    Let $\M\subset M$ denote the subset corresponding to elastic collisions between the two balls. It is shown in \cite[Proposition 1.8]{BBNV} that, with the first return map to $\M$, this system can be modeled with a two-sided CMZ structure, furthermore, our condition (C1) applies. The argument of \cite{BBNV} relies on unstable curves and associated standard pairs, in the sense of (C2). It is also proved that $\mu_{\dhat}(H_n)\ll n^{-3}$.
    Moreover, in \cite[Sections 2.3-2.5]{BBNV} the authors introduce a simplified model and show that it applies to the system. Using this, a lower bound $\mu_{\dhat}(H_n)\gg n^{-3}$ can be established using the same argument and calculations as in \cite[Section 2.3]{BBNV}.
    Nonetheless, \cite{BBNV} proves a weaker tail bound for the return times to the base of the Young tower, in particular
    \[
        \sum_{k\geq n}\mu_{\dnull}(A_k)\ll \frac{(\log n)^3}{n^2}.
    \]
    On the other hand, in \cite[Classification 4.2]{BBNV} it is concluded that $R_k$ has a first step expansion rate $k^2$, and the authors prove a Growth Lemma in \cite[Section 4]{BBNV}, which means that our condition (C5) is satisfied with $p=2$.
    The calculations in \cite[Section 2.3]{BBNV} about the measure of $\mu(R_k)$ show that $t=3$, establishing (C3) and (C4).

    This means that all (C1)-(C5) hold, thus by Theorem \ref{main},
    \[
        \mu_{\dnull}(A_n)\asymp n^{-3}.
    \]

    Consequently, by Theorem \ref{correlation general},
        \[
            |C_T(f, g, n)|
            \ll
            n^{-2}
        \]
    for all $f, g\in \mathcal{H}(M)$.

     By Theorem \ref{correlation},
        \[
            |C_T(f, g, n)|
            \asymp
            n^{-2}
        \]
    for all $f, g\in \mathcal{H}(\M)$ with $\int_M f \operatorname{d}\! \mu\int_M g \operatorname{d}\!\mu \neq 0$ and
        \[
            |C_T(f, g, n)|
            \ll
             n^{-3}\log n
        \]
    for all $f, g\in \mathcal{H}(\M)$ with $\int_M f \operatorname{d}\! \mu=0$.

    Furthermore, we can improve the former ASIP rate estimate $o(n^{1/3}(\log n)^{4/3 + \varepsilon})$ of \cite[Theorem 1.1]{ASIP} to a rate $o(n^{1/3}(\log n)^{1/3 + \varepsilon})$ (for any $\varepsilon>0)$.
\end{example}

\begin{example}{\textbf{Bunimovich flowers}}\label{flowers}

Billiard models study the motion of a point particle in some domain
$Q\subset \mathbb{R}^2$. The particle has unit speed, within $Q$
the motion is uniform, while reflections at the boundary
$\partial Q$ are specular (angle of reflection equals angle
of incidence).
Let $M=\partial Q\times \left[-\frac{\pi}{2},\frac{\pi}{2}\right]$, the phase space corresponding to collisions at the boundary, where the second coordinate measures the angle the outgoing velocity makes with the normal vector of the boundary at the point of reflection. This way the billiard map $T:M\to M$ arises for which there is a natural invariant probability measure absolutely continuous with respect to Lebesgue measure, the Liouville measure. We refer to \cite{billiards} for further details on billiards.
Specifically, let us consider Buminovich-type billiard tables as in
\cite[Section 8]{ChZh} or \cite[Example 8.9]{Ian}.
This means that the boundary $\partial Q=\Gamma=\Gamma_1\cup\dots\cup\Gamma_r$  is a simple, closed, piecewise $C^3$-smooth
curve such that each smooth component $\Gamma_i$ has nonvanishing curvature, and there is at least one convex inwards (dispersing) component, while the convex outwards (focusing) components are circular arcs strictly smaller than a semicircle.
Furthermore, we assume that all corner points have positive angle,
and each convex outwards arc can be continued to a circle which is contained entirely inside $Q$.
These billiard tables were introduced by Bunimovich in \cite{Bu} who also established that the corresponding billiard map $T:M\to M$ is hyperbolic and ergodic with respect to the Liouville measure.

To proceed, following \cite[Section 8]{ChZh} let us introduce $\M\subset M$ which consist of (i) all the collisions at dispersing boundary components, and (ii) only the first collisions at focusing boundary components. It is proved in \cite[Theorem 10]{ChZh} that the first return map to $\M$ enjoys exponential decay of correlations and $T:M\to M$ possesses a two-sided Chernov-Markarian-Zhang structure, thus (C1) holds.

Furthermore, in \cite[Section 8]{ChZh} the authors actually check that all our further conditions (C2)-(C5) hold for these billiard systems. Indeed, the method of \cite{ChZh} relies on the use of standard pairs, thus one can see that  (C2) and (C3) are satisfied.
\cite[Section 8]{ChZh} also checks that unstable sizes of the sets $R_k$ can be dominated by $k^{-2}$, thus (C4) holds.
In \cite[Section 8]{ChZh} and \cite[Chapter 8]{billiards} it is also established that the first step expansion rate on unstable curves in $R_k$ is at least $c\cdot k^{3/2}$ for some $c>0$\footnote{The expansion rate is $k^2$ in the p-metric, but only $k^{3/2}$ can be guaranteed in the Euclidean metric which is relevant for us.} \cite[Section 8]{ChZh} also proves
the classical version of the Growth Lemma (ie.~Lemma~\ref{l:growthclassic}) for this class of Bunimovich billiards. Thus we conclude  that (C5) holds with $p=3/2$ and $t=2$.
On \cite[page 30]{ChZh} it is shown that $\mu_{\dhat}(H_n) \ll n^{-3}$.
\cite[Example 8.9]{Ian} includes lower bounds on $\mu(H_n)$ of the same order, thus we have
\[
    \mu_{\dhat}(H_n)\asymp n^{-3}.
\]
This altogether means that we can apply Theorem \ref{main} to show that
\[
    \mu_{\dnull}(A_n)\asymp n^{-3}.
\]

    By Theorem \ref{correlation general},
        \[
            |C_T(f, g, n)|
            \ll
            n^{-2}
        \]
    for all $f, g\in \mathcal{H}(M)$.

 By Theorem \ref{correlation},
    \[
    |C_T(f, g, n)|
    \asymp
    n^{-2}
    \]
    for all $f, g\in \mathcal{H}(\M)$ with $\int_M f \operatorname{d}\! \mu\int_M g \operatorname{d}\!\mu \neq 0$ and
    \[
        |C_T(f, g, n)|
        \ll
         n^{-3}\log n
    \]
    for all $f, g\in \mathcal{H}(\M)$ with $\int_M f \operatorname{d}\! \mu=0$.

    We can also improve the former ASIP rate estimate $o(n^{1/3}(\log n)^{4/3 + \varepsilon})$ of  \cite[Theorem 1.1]{ASIP} to a rate of $o(n^{1/3}(\log n)^{1/3 + \varepsilon})$, for any $\varepsilon>0$.

\end{example}

\begin{example}{\textbf{Billiards with flat points}}\label{flat points}

    In \cite{variable mixing rates} the authors consider a family of planar billiards that depend on a parameter $\beta>2$.
    The boundary is given by the functions $y=x^{\beta} \pm 1$ in the neighbourhood of the points $(0, \pm 1)$. Hence the billiard map has a period 2 orbit bouncing perpendicularly back and forth between two points of the boundary with vanishing curvature.
    Let $\M\subset M$ be the set of collisions that occur outside a small neighbourhood of these two points.
    By \cite[property \textbf{(F1)}]{variable mixing rates}, the system with the first return map of $\M$ satisfies (C1). Although \cite{variable mixing rates} focuses on the evolution of unstable manifolds, it is clear from the exposition that the analysis applies to arbitrary unstable curves, which means that properties (C2) and (C3) apply. By \cite[property \textbf{(F2)}]{variable mixing rates},
    $\mu_{\dhat}(H_n) \ll n^{-\alpha}$ with $\alpha=\frac{\beta + 2}{\beta - 2}+1$.
    \cite[inequality (3.3)]{variable mixing rates} is essentially a Growth Lemma in the sense of Lemma~\ref{l:growthclassic}.
    Furthermore, the first step expansion rate is asymptotic to $n^{\alpha+1}$, while the unstable sizes of $R_n$ can be dominated by $n^{-\alpha-1}$ according to \cite[Proposition 2]{variable mixing rates}, hence (C5) and (C4) apply.
    Below the proposition, we can see that
    \[
        \mu_{\dhat}(H_n)\asymp n^{-\alpha} = n^{-\frac{\beta + 2}{\beta - 2}-1}.
    \]
    Hence, by Theorem~\ref{main},
    \[
        \mu_{\dnull}(A_n)\asymp n^{-\alpha} = n^{-\frac{\beta + 2}{\beta - 2}-1}.
    \]

    By Theorem \ref{correlation general},
        \[
            |C_T(f, g, n)|
            \ll n^{-\alpha+1}=
            n^{-\frac{\beta + 2}{\beta - 2}}
        \]
    for all $f, g\in \mathcal{H}(M)$.

     By Theorem \ref{correlation},
        \[
        |C_T(f, g, n)|
        \asymp
        n^{-\frac{\beta + 2}{\beta - 2}}
        \]
    for all $f, g\in \mathcal{H}(\M)$ with $\int_M f \operatorname{d}\! \mu\int_M g \operatorname{d}\!\mu \neq 0$ and
        \[
            |C_T(f, g, n)|
            \ll
             n^{-\frac{\beta + 2}{\beta - 2}-1}\log n
        \]
    for all $f, g\in \mathcal{H}(\M)$ with $\int_M f \operatorname{d}\! \mu=0$.

        See \cite{Bonnafoux} for similar results on the associated billiard flow.

        Concerning ASIP rates for the billiard map, let us first mention \cite{Chen} which studies the question with a different technique, and in particular the result of \cite[Section 4.2]{Chen} on ASIP for our Example~\ref{flat points} with rate $o(n^{\lambda})$ for any $\lambda>\max\{1/4,1/\alpha\}$. Although not  explicitly stated, the methods of \cite{ASIP} also apply to this example (as mentioned in the abstract of that paper) and already improve upon the result of \cite{Chen}. Nonetheless, by Theorem~\ref{ASIP theorem} we get a further improved rate $o(n^{1/\alpha}(\log n)^{1/\alpha + \varepsilon})$ for the billiard map, for arbitrary small $\varepsilon>0$ (recall $\alpha>2$).

\end{example}

\appendix

\section{Appendix}

\subsection{Young towers and Chernov-Markarian-Zhang structures}\label{CMZ structures}
Here we describe the framework we use in the paper, namely, two-sided Chernov-Markarian-Zhang structures. We include the necessary notions needed to define these structures. For full details, see \cite{Ian}. In the following table we summarize the notation differences between our paper and \cite{Ian}.
\\[1cm]
\begin{center}
\begin{tabular}
{>{\arraybackslash}m{4.5cm}|>{\centering\arraybackslash}m{3.5cm}|>{\centering\arraybackslash}m{3.5cm}}
\hline
    \textbf{Meaning} & \textbf{Our notation} & \textbf{Notation in \cite{Ian}} \\
    \hline
    The system                                  & $T: M \to M$                  & $f: M\to M$ \\
    The subset                                  & $\M \subset M$                & $X \subset M$ \\
    It's first return time                      & $R : \M \to \N$             & $h: X \to \N$\\
    The subsystem                               & $\T=T^R:\M \to \M$            & $f_{X}=f^h: X \to X$ \\
    The base of the tower                       & $\dnull \subset \M$           & $Y \subset X$ \\
    It's return time in the subsystem           & $\sigma: \dnull \to \N$     & $\sigma: Y \to \N$ \\
    The return map of the base                  & $F_0: \dnull \to \dnull$      & $F: Y \to Y$ \\
    The partition on the base                   & $\mathcal{P}$                 & $\alpha$\\
    The inner tower                             & $\dhat=\dnull^{\sigma}$       & $\Delta_{rapid}=Y^{\sigma}$ \\
    The inner tower map                         & $\F: \dhat \to \dhat$         & no notation \\
    The return time of the base in the system   & $h: \dnull \to \N$          & $\varphi: Y \to \N$ \\
    The whole tower                             & $\Delta = \dnull^h$           & $ \Delta=Y^{\varphi}$ \\
    The tower map                               & $F: \Delta \to \Delta$        & $f_{\Delta}: \Delta \to \Delta$ \\
    The one-sided version of the tower map (in \ref{2sided CMZ})     & $F': \Delta' \to \Delta'$     & $\bar{f}_{\Delta}: \bar{\Delta} \to \bar{\Delta}$\\
    \hline
\end{tabular}
\\[1cm]
\end{center}
First we describe the one-sided versions of these notions and then extend them to invertible systems as well.
\subsubsection*{Gibbs-Markov Maps}
Suppose that $(\dnull, \mu_{\dnull})$ is a probability space with a countable measurable partition $\mathcal{P}$, and let $F_0: \dnull \to \dnull$ be an ergodic measure-preserving transformation.
Define the separation time $s(x, y)$ to be
\[
    s(x, y)=\min\set{n\in\N}{\exists A\in \mathcal{P}:\; F_0^n(x)\in A, F_0^n(y)\notin A }.
\]
It is assumed that the partition $\mathcal{P}$ separates trajectories, so $s(x, y) = \infty$ if and only if $x = y$.
Then $\theta^s$ is a metric for any $\theta \in (0, 1)$.
\\
Let
\[
\xi = \frac{d\mu_{\dnull}}{d\mu_{\dnull} \circ F_0}: \dnull \to \mathbb{R}.
\]
We say that $F_0$ is a (full-branch) Gibbs-Markov map if:
\begin{itemize}
    \item $F_0|_A: A \to \dnull$ is a measurable bijection for each $A \in \mathcal{P}$, and
    \item There are constants $C > 0$, $\theta \in (0, 1)$ such that
    \[
    | \log \xi(x) - \log \xi(y) | \leq C \theta^{s(x, y)}
    \]
    for all $x, y \in A$, $A \in \mathcal{P}$.
\end{itemize}
Hence, there is a constant $C > 0$ such that
\[
\xi(x) \leq C \mu_{\dnull}(A) \quad \text{and} \quad | \xi(x) - \xi(y) | \leq C \mu_{\dnull}(A) \theta^{s(x, y)}, \tag{2.1}
\]
for all $x, y \in A$, $A \in \mathcal{P}$.
\subsubsection*{Return Maps}
Let $(M, \mu)$ be a probability space and $T: M \to M$ an ergodic measure-preserving transformation.
Let $\M \subset M$ be a measurable subset with $\mu(\M) > 0$. If $R: \M \to \N$ is integrable and $T^{R(x)}x \in \M$ for all $x \in \M$, then $R$ is called a \emph{return time} and $T^R: \M \to \M$ is called a \emph{return map}.
\\
If $R$ is the first return time of $\M$ under $T$ (i.e., $R(x) = \inf\{n \geq 1 \,| \, T^n x \in \M\}$), then $T^R: \M \to \M$ is called the \emph{first return map} and
\[
\mu_{\M} = \frac{\mu|_{\M}}{\mu(\M)}
\]
is an ergodic $T^R$-invariant probability measure on $\M$.
\\
\subsubsection*{Young Towers}
Let $F_0: \dnull \to \dnull$ be a full-branch Gibbs-Markov map on $(\dnull, \mu_{\dnull})$ with partition $\mathcal{P}$ and let $h: \dnull \to \N$ be an integrable function constant on partition elements. The (one-sided) Young tower $\Delta = \dnull^h$ and tower map $F: \Delta \to \Delta$ can be defined the following way:
\[
\Delta = \{(x, \ell) \in \dnull \times \N : 0 \leq \ell < h(x)\}, \quad
F(x, \ell) =
\begin{cases}
(x, \ell + 1) & \text{if } \ell < h(x) - 1, \\
(F_0x, 0) & \text{if } \ell = h(x) - 1.
\end{cases}
\]
Let $\bar{h} = \int_{\dnull} h \, d\mu_{\dnull}$ and let $\nu$ be the counting measure on $\N$.
Then there is an ergodic $F$-invariant probability measure
\[
\mu_\Delta = \frac{\mu_{\dnull} \times \nu }{\bar{h}}
\]
 on $\Delta$.
\\
Let $T: M \to M$ be an ergodic measure-preserving transformation on a probability space $(M, \mu)$, and let $\dnull \subset M$ be measurable with $\mu(\dnull) > 0$.
Suppose that  $h: \dnull \to \N$ is a return time, constant on partition elements, and that $F_0=T^h: \dnull \to \dnull$ is a full-branch Gibbs-Markov map with respect to a probability measure $\mu_{\dnull}$ on $\dnull$.
Construct the tower $\Delta = \dnull^h$ and tower map $F: \Delta \to \Delta$.
\\
If the semiconjugacy $\pi_M: \Delta \to M$, $\pi_M(y, \ell) = T^\ell y$ is measure preserving as well,
then we say that $T$ is \emph{modelled by a Young tower}.
\\
\subsubsection*{Chernov-Markarian-Zhang Structure}
Suppose $\dnull \subset \M \subset M$ are Borel sets with $\mu(\dnull) > 0$.
Denote by $R: \M \to \N$ and  $\T = T^{R}: \M \to \M$ the first return time and first return map of $\M$.
\\
We assume that $\T: \M \to \M$ is modelled by a Young tower $\dhat = \dnull^\sigma$ with return time $\sigma: \dnull \to \N$ and return map $F_0 = \T^\sigma: \dnull \to \dnull$ which is a full-branch Gibbs-Markov map .
We require  $R$ to be constant on $\T^{\ell} A$ for all $A \in \mathcal{P}$, $0 \leq \ell < \sigma(A)$, thus there is an integrable return time $h$ such that $T: M \to M$ is modelled by a Young tower $\Delta = \dnull^h$ with the same Gibbs-Markov map $F_0 = \T^{\sigma} = T^h$.

If $T: M \to M$ satisfies these assumptions, we say it has a \emph{Chernov-Markarian-Zhang structure}.

\subsection{Two-sided Young towers and Chernov-Markarian-Zhang structures}\label{2sided CMZ}
Now we introduce a two-sided (invertible) analogue of the structures described previously.
Throughout the Section, $T : M \to M$, $\T : \M \to \M$, $F : \Delta \to \Delta$ and $F_0 : \dnull \to \dnull$ are the two-sided versions of the maps from Section \ref{CMZ structures}, and the one-sided versions are denoted $F_0' : \dnull' \to \dnull'$ and so on.
\subsection*{Two-sided Gibbs-Markov Maps}
Suppose that $(\dnull, d)$ is a bounded metric space with a Borel probability measure $\mu_{\dnull}$ and let $F_0 : \dnull \to \dnull$ be an ergodic measure-preserving transformation. Let $F_0' : \dnull' \to \dnull'$ be a full-branch Gibbs-Markov map with partition $\mathcal{P}$ and ergodic invariant probability measure $\mu'_{\dnull}$.
\\
We suppose that there is a measure-preserving semiconjugacy $\pi' : \dnull \to \dnull'$, which makes the following diagram commutative.

\[
\begin{tikzcd}
\dnull \arrow[r, "F_0"] \arrow[d, "\pi'"'] &\dnull \arrow[d, "\pi'"] \\
\dnull' \arrow[r, "F_0'"'] & \dnull'
\end{tikzcd}
\]

Then the separation time on $\dnull'$ lifts to a separation time on $\dnull$ defined by $s(x, y) = s(\pi'x, \pi'y)$ for $x, y \in \dnull$. We assume that there are constants $C > 0$, $\theta \in (0, 1)$ such that
\begin{equation}\label{eq:whw}
d(F_0^n x, F_0^n y) \leq C(\theta^n + \theta^{s(x,y)-n}) \quad \text{for all } x, y \in \dnull, \, n \geq 1.
\end{equation}
If $F_0 : \dnull \to \dnull$ satisfies all the assumptions above, then it is called a \textit{two-sided Gibbs-Markov map}.
\begin{remark}\label{r:whw}
Note that instead of \eqref{eq:whw} it would be also very natural to assume a slightly stronger inequality, namely
\begin{equation}\label{eq:whw1}
d(F_0^n x, F_0^n y) \leq C(\theta^n d(x,y) + \theta^{s(x,y)-n}) \quad \text{for all } x, y \in \dnull, \, n \geq 1.
\end{equation}
The advantage would be that \eqref{eq:DistIneqAlexey} together with \eqref{eq:whw1} implies \eqref{eq:DistIneqIan}.
\end{remark}

\subsection*{Two-sided Young Towers}
Suppose that $F_0 : \dnull \to \dnull$ is a two-sided Gibbs-Markov map on $(\dnull, \mu_{\dnull})$ and let $h : \dnull \to \N$ be an integrable function which is constant on $\pi'^{-1}A$ for each $A \in \mathcal{P}$.
This way, $h$ is well-defined on $\dnull'$.
Define the one-sided Young tower $\Delta' = \dnull'^h$ and tower map $F' : \Delta' \to \Delta'$ as in Section \ref{CMZ structures}.
Using $F_0 : \dnull \to \dnull$ instead of $F_0' : \dnull' \to \dnull'$, we also define the two-sided Young tower $\Delta = \dnull^h$ and tower map $F : \Delta \to \Delta$.
There are ergodic invariant probability measures
\[
\mu_\Delta = \frac{\mu_{\dnull} \times \nu}{ \bar{h}}, \quad \mu'_\Delta = \frac{\mu'_{\dnull} \times \nu}{ \bar{h}},
\]
on $\Delta$ and $\Delta'$. Recall, that $\nu$ is the counting measure on $\N$.
\\
We can extend the projection $\pi' : \dnull \to \dnull'$ to $\pi' : \Delta \to \Delta'$ with $\pi'(y, \ell) = (\pi'y, \ell)$. Then $\pi'$ is a measure-preserving semiconjugacy between $F$ and $F'$.
\\
Now assume that $T : M \to M$ is an ergodic measure-preserving transformation on a probability space $(M, \mu)$, and that $\dnull \subset M$ is measurable with $\mu(\dnull) > 0$. Suppose that $F_0 : \dnull \to \dnull$ is a two-sided Gibbs-Markov map with respect to a probability measure $\mu_{\dnull}$ on $\dnull$, and that $h : \dnull \to \N$ is a return time as above.
Construct the tower $\Delta = \dnull^h$ and the tower map $F : \Delta \to \Delta$.
If the semiconjugacy $\pi_M : \Delta \to M$, $\pi_M(y, \ell) = T^\ell y$ is also measure preserving, then we say that $T$ is modelled by a \textit{two-sided Young tower}.
\subsection*{Two-sided Chernov-Markarian-Zhang Structure}
Suppose that $(M, d)$ is a bounded metric space with Borel probability measure $\mu$ and let $T : M \to M$ be a mixing measure-preserving transformation.
Suppose that $\dnull \subset \M \subset M$ are Borel sets satisfying $\mu(\dnull) > 0$.
Denote the first return time by $R : \M \to \N$ and first return map by $\T = T^R : \M \to \M$.
\\
We require that $\T : \M \to \M$ is modelled by a two-sided Young tower $\dhat = \dnull^\sigma$ with return time $\sigma : \dnull \to \N$ and return map $F_0 = \T^\sigma : \dnull \to \dnull$, which is a two-sided Gibbs-Markov map with ergodic invariant probability measure $\mu_{\dnull}$ and partition $\mathcal{P}$ such that $\sigma$ is constant on partition elements.
We further assume that $R$ is constant on $\T^\ell \pi'^{-1}A$ for all $A \in \mathcal{P}$, $0 \leq \ell < \sigma(A)$.
\\
Then there is an integrable return time $h = R_\sigma : \dnull \to \N$ such that $T : M \to M$ is modelled by a Young tower $\Delta = \dnull^h$ with the same two-sided Gibbs-Markov map $F_0 = \T^\sigma = T^h$.
\\
If $T : M \to M$ satisfies these assumptions, it is said to have a \textit{two-sided Chernov-Markarian-Zhang structure}.

\subsection{Proof of some properties of regularly varying functions}\label{proof of rv}
\begin{customproof}{Proof of Lemma \ref{index}}
    We will use Karamata's integral representation theorem of regularly varying functions, see \cite[Theorem 1.3.1]{regular variation} for example.
\begin{theorem}
    A function $r:\mathbb{R}^+ \rightarrow \mathbb{R}^+$ is regularly varying, if and only if there exist measurable functions $c, a:\mathbb{R}^+ \rightarrow \mathbb{R}$ with $\exists C=\lim_{x\rightarrow \infty}c(x) \in \mathbb{R}$, $\exists -\alpha =  \lim_{x\rightarrow \infty}a(x)\in\mathbb{R} $ and there is a $y>0$ such that
    \[
        r(x) = \exp\left(c(x)+\int_y^{x}\frac{a(s)}{s}\operatorname{d}s\right)
    \]
    holds for all $x\geq y$.
\end{theorem}

With this, we can see that for any $\varepsilon>0$ there exists an $x_0>0$ such that for all $x>x_0$,
\begin{align*}
    \log r(x)
    \leq&
    C+\varepsilon + \int_y^{x}\frac{a(s)}{s}\operatorname{d}s
    \leq
    C+\varepsilon + \int_y^{x_0}\frac{a(s)}{s}\operatorname{d}s + \int_{x_0}^{x}\frac{a(s)}{s}\operatorname{d}s
    \\
    \leq&
    C+\varepsilon + \int_y^{x_0}\frac{a(s)}{s}\operatorname{d}s + (-\alpha+ \varepsilon)\int_{x_0}^{x}\frac{1}{s}\operatorname{d}s
    \\
    =&
    C+\varepsilon + \int_y^{x_0}\frac{a(s)}{s}\operatorname{d}s + (-\alpha+ \varepsilon)(\log x - \log x_0)
    \\
    =&
    C+\varepsilon + \int_y^{x_0}\frac{a(s)}{s}\operatorname{d}s - (-\alpha+ \varepsilon) \log x_0 +(-\alpha+ \varepsilon)\log x
\end{align*}
Hence,
\begin{align*}
    \frac{\log r(x)}{\log x}
    \leq
    -\alpha+\varepsilon
    +
    \left( C+\varepsilon + \int_y^{x_0}\frac{a(s)}{s}\operatorname{d}s - (-\alpha+ \varepsilon) \log x_0 \right)(\log x)^{-1}
\end{align*}
Which means that there is an $x_1>0$ such that for all $x>x_1$,
\[
    \frac{\log r(x)}{\log x}
    \leq
    -\alpha + 2 \varepsilon
\]
holds. By using similar calculations, one can arrive at the lower bound
\[
    -\alpha - 2 \varepsilon
    \leq
    \frac{\log r(x)}{\log x}
\]
From this, we can see that
\[
    \exists \lim_{x \rightarrow \infty} \frac{\log r(x)}{\log x} =  -\alpha.
\]
We can see that $-\alpha$ is the index of $r$
\end{customproof}

\begin{customproof}{Proof of Lemma \ref{asymptotic}}
Since $r$ is monotone decreasing,
    \[
    \sum_{k\geq n} r(k)
    \geq
    \sum_{k\geq n}\int_k^{k+1} r(x) \operatorname{d} x
    =\int_n^{\infty}r(x) \operatorname{d} x.
    \]
    On the other hand,
    \begin{align*}
        \sum_{k\geq n} r(k)
        =&
        r(n)+ \sum_{k\geq n}\int_k^{k+1} r(k+1) \operatorname{d} x
        \\
        \leq&
        r(n)+ \sum_{k\geq n}\int_k^{k+1} r(x) \operatorname{d} x
        =
        \int_n^{\infty}r(x) \operatorname{d} x + r(n)
        \\
        \ll&
        \int_n^{\infty}r(x) \operatorname{d} x + n^{-\alpha}
        \ll
        \int_n^{\infty}r(x) \operatorname{d} x + n^{-\alpha-\varepsilon +1}
        \\
        \ll&
        \int_n^{\infty}r(x) \operatorname{d} x
        ,
    \end{align*}
where we have used that $\int_n^{\infty}r(x) \operatorname{d} x\gg n^{-\alpha-\varepsilon +1}$ by
\eqref{almost poly}.

\end{customproof}

\end{document}